\documentclass[12pt]{article}
\usepackage{amsmath}
\usepackage{mathtools}
\usepackage{amssymb}

\usepackage{bm} 
\usepackage{booktabs} 

\usepackage{hyperref} 

\newtheorem{remark}{Remark}

\hypersetup{pdfauthor={Roy R. Lederman},pdftitle={Generalized Prolate Spheroidal Functions}}

\title{Numerical Algorithms for the Computation of Generalized Prolate Spheroidal Functions}
\date{\today}
\author{Roy R. Lederman\\ \small The Program in Applied and Computational Mathematics\\ \small Princeton University}

\begin{document}
\maketitle

\begin{abstract}
  Generalized Prolate Spheroidal Functions (GPSF) are the eigenfunctions of the truncated Fourier transform, restricted to D-dimensional balls in the spatial domain and frequency domain. Despite their useful properties in many applications, GPSFs are often replaced by crude approximations. 
  The purpose of this paper is to review the elements of computing GPSFs and associated eigenvalues.
  This paper is accompanied by open-source code.   
  \end{abstract}

%
%

\section{Introduction}\label{sec:intro}

Prolate Spheroidal Wave Functions (PSWF), the one-dimensional case of the GPSFs discussed in this paper, 
are the eigenfunctions of the truncated Fourier transform $F_c$, defined by the formula
\begin{equation}
  \left(F_c\psi\right)(x)   = \int_{-1}^1 e^{\mathrm{i} cxy} \psi(y) dy ,
\end{equation}
with $x \in [-1,1]$. 
In the classic work \cite{slepian1961prolate} by Slepian and Pollak it has been shown that the integral operator $F_c ^*F_c$ 
commutes with the differential operator $L_c$ defined by the formula
\begin{equation}
   \left(L_c\psi\right)(x) = (1-x^2) \frac{d^2}{dx^2}\psi(x) - 2x \frac{d}{dx}\psi(x) -c^2 x^2 \psi(x).
\end{equation}
This remarkable property was subsequently used to infer many of the properties of PSWFs and to develop efficient and accurate methods for computing them in \cite{slepian1961prolate,xiao2001prolate,xiao2003high,rokhlin2007approximate,osipov2013prolate,lederman2017lower} among other works. 
This rare relation between an integral operator and a differential operator has been discovered and used for the analysis of the truncated Laplace transform as well (e.g. \cite{bertero1985commuting,LedermanLaplace2,lederman2017lower}).

In subsequent work in \cite{slepian1964prolate}, Slepian discovered a similar relation associated with the general case of 
high dimensional Fourier transform, supported on the unit ball
\begin{equation}
 \left( {F}_c \psi\right)({\bm x}) =    \int_{|{\bm y}| \leq 1} e^{ \mathrm{i} c \left( {\bm{x}} \cdot {\bm{y}} \right) } \psi(\bm{y})  {d}\bm{y} ,
\end{equation}
where ${\bm y}$ and ${\bm x}$ are vectors in $\mathbb{R}^D$, and $ |{\bm x}| \leq 1 $ .
The eigenfunctions of this operator are the GPSFs. 

PSWFs have been useful in many applications in signal processing. Among their many useful properties, they are the optimal basis for representing functions whose energy is as concentrated as possible in the interval $[-1,1]$ in the spatial domain, and in the interval $[-c,c]$ in the frequency domain. 
GPSFs are known to be useful in many applications in signal processing and optics, however they have gained a reputation of being prohibitively difficult to compute, and in practice they are often replaced by crude approximations. 

The properties of GPSFs have been studied by Slepian in \cite{slepian1964prolate}. Many of the properties of the two-dimensional GPSFs, have been studied in \cite{shkolnisky2007prolate,shkolnisky2006approximation}.
Recently, the more general GPSFs have been studied in \cite{serkh2015generalized}, and they are being studied in the upcoming \cite{greengard2017generalized}.
However, we are not aware of a complete publicly available resource and code for efficient and accurate computation of GPSFs and their associated eigenvalues.

The purpose of this paper is to review the essential results for computing GPSFs, and to introduce numerical code for computing GPSFs.
We review results primarily from \cite{slepian1964prolate,shkolnisky2007prolate} and \cite{serkh2015generalized}, with some reformulation and generalization, as well as some unpublished results which Philip Greengard and Kirill Serkh have been generous to share with us from their upcoming \cite{greengard2017generalized}, and introduce additional  results.

The immediate motivation for this paper is the use of three-dimensional GPSFs for representing objects that are highly concentrated in both the spatial and frequency domains in cryo-electron microscopy (cryo-EM) applications (e.g. \cite{lederman2016representation,lederman2017continuously}). 
Indeed the two-dimensional GPSFs have been used successfully for the representation of images in cryo-EM in \cite{landa2017approximation,landa2017steerable}. 

A MATLAB\textsuperscript{TM} source code is available at \url{http://github.com/lederman/prol} together with the LaTeX source for this paper. FORTRAN code is under development.
In addition, the URL above contains Mathematica\textsuperscript{TM} code which performs various levels of analytic verification of some of the equations in this paper after they have been derived in other means (an ``open-source proof''). 

The current paper and code are preliminary results of the ongoing project in  \url{http://github.com/lederman/prol}. We welcome comments and contributions to the code and paper. 

This paper is organized as follows. Section \ref{sec:pre} briefly reviews standard definitions and results that are used in the paper. Section \ref{sec:analysis} presents the essential analytical results required for computing GPSFs. 
The algorithms based on these analytical results are presented in Section \ref{sec:alg}.
Numerical results obtained with the publicly available code are presented in Section \ref{sec:results}.
Finally, brief conclusions are presented in \ref{sec:conclusions}.

\section{Preliminaries}\label{sec:pre}

%
%
%

\begin{table}[htbp]
  \caption{Notation}
\begin{center}
\begin{tabular}{r p{10cm} }
\toprule
$c$ &                           band-limit, see, for example, Equation (\ref{eq:def:truncated_Fourier}). \\
$D$ &                           dimensionality. $D \ge 2$ is assumed.\\
$p$ &                           $p=D-2$ \\
$N$ &                           angular frequency, see, for example, Section \ref{sec:surfaceharmonics} and equations (\ref{eq:scaled_truncated_fourier_eig:separation}) and (\ref{eq:def:prol}) \\
$P^{(\alpha,\beta)}_n$ &         Jacobi Polynomials (see Section \ref{sec:jacobi}  )\\
$\overline{R_{\text{N,n}}^{\text{p}}}$ & normalized radial Zernike polynomials (see Section \ref{sec:zernike}) \\
$T_{\text{N,n}}^{\text{p}}$ &     weighted radial Zernike polynomials (see Section \ref{sec:zernike}) \\
$J_n$ &                       Bessel functions of the first kind (see Section \ref{sec:bessel}) \\
$S_N^m({\bm \xi})$ &            surface harmonics (see Section \ref{sec:surfaceharmonics}) \\
$\Gamma (n)$  &                 the gamma function. Defined in (\ref{eq:def:gammafunc})\\
$\mathcal{F}$ &                 Fourier transform. Defined in (\ref{eq:def:Fourier}) \\
$\mathcal{F}_c$ &               truncated Fourier transform. Defined in (\ref{eq:def:truncated_Fourier}) \\
${F}_c$ &                       scaled truncated Fourier transform. Defined in (\ref{eq:scaled_truncated_fourier}) \\
$M_{p,c,N}$ &                    the integral operator defined in (\ref{eq:operator_M}) \\
$L_{p,c,N}$ &                    the differential operator defined in (\ref{eq:operator_L}) \\
$B_{\text{N}}^{\text{p,c}}$ &       the matrix form of the operator $L_{p,c,N}$ in the basis of weighted radial Zernike polynomials. Defined in (\ref{eq:Bmat_nn}) and (\ref{eq:Bmat_nm1n})\\
$\psi_{N,n,m}({\bm r})$ &        GPSFs defined in (\ref{eq:def:prol}) \\
$\Phi_{N,n}(r)$ &               radial GPSFs defined in (\ref{eq:def:radial_prol}) \\
$\varphi_{N,n}(r)$ &            weighted radial GPSF defined in (\ref{eq:def:weighted_radial_prol}) \\
$\widetilde{\psi_{N,n,m}}({\bm r})$ &    the ``eigenfunctions'' defined in (\ref{eq:def:tilde_psi}) \\
$\alpha_{N,n,m} = \alpha_{N,n}$ & the eigenvalues defined in (\ref{eq:alpha}) \\
$\beta_{N,n}$ &                 the eigenvalues defined in (\ref{eq:def:radial_prol}) \\
$\gamma_{N,n}$ &                the eigenvalues defined in (\ref{eq:gamma}) \\
$\nu_{N,n}$ &                   the eigenvalues defined in (\ref{eq:nu}) \\
$\chi^{p,c,N}_n$ &               the eigenvalues defined in (\ref{eq:operator_L:eig}) of the differential operator $L_{p,c,N}$  \\
\bottomrule
\end{tabular}
\end{center}
\label{table;notation}
\end{table}

\subsection{Jacobi Polynomials}\label{sec:jacobi}
The following are well-known properties of the Jacobi polynomials, denoted by $P_n^{(\alpha ,\beta )}(x)$. 
The standard definition of Jacobi polynomials, along with these properties and others, can be found, inter alia, in \cite{abramowitz1964handbook}.

\begin{equation}\label{eq:jac1}
(1-x) P_n^{(\alpha +1,\beta )}(x)=\frac{(n+\alpha +1) P_n^{(\alpha ,\beta )}(x)-(n+1)
   P_{n+1}^{(\alpha ,\beta )}(x)}{n+\frac{\alpha }{2}+\frac{\beta }{2}+1}
\end{equation}

\begin{equation}\label{eq:jac2}
(n+\alpha +\beta ) P_n^{(\alpha ,\beta )}(x)=(2 n+\alpha +\beta ) P_n^{(\alpha ,\beta
   -1)}(x)-(n+\alpha ) P_{n-1}^{(\alpha ,\beta )}(x)
\end{equation}

\begin{equation}\label{eq:jac3}
\frac{{d} P_n^{(\alpha ,\beta )}(x)}{{d} x}=\frac{1}{2} (n+\alpha +\beta +1)
   P_{n-1}^{(\alpha +1,\beta +1)}(x)
\end{equation}

\begin{equation}\label{eq:jac4}
  P_n^{(\alpha ,\beta )}(1) = \frac{\Gamma(n+\alpha+1)}{\Gamma(n+1)\Gamma(\alpha+1)}
\end{equation}

\subsection{Radial Zernike Polynomials}\label{sec:zernike}

The purpose of this section is to define the normalized radial Zernike polynomials.
The properties of Zernike polynomials are discussed in further detail in \cite{slepian1964prolate,serkh2015generalized}. 
We note that slightly different definitions and normalization are used in different sources. 

In this paper, the normalized radial Zernike polynomials, denoted by $\overline{R_{\text{N,n}}^{\text{p}}}(x)$, are defined by the formula
\begin{equation}\label{eq:def:norm_radial_zer}
\overline{R_{\text{N,n}}^{\text{p}}}(x)=\sqrt{2} (-1)^n x^N \sqrt{2 n+N+\frac{p}{2}+1}
   P_n^{\left(N+\frac{p}{2},0\right)}\left(1-2 x^2\right),
  \end{equation}
where $P_n^{(\alpha ,\beta )}(x)$ are the Jacobi polynomials. 
The normalized radial Zernike polynomials are orthonormal in the following sense: 
\begin{equation}
\int_0^1 x^{p+1} \overline{R_{\text{N,n}}^{\text{p}}}(x) \overline{R_{\text{N,j}}^{\text{p}}}(x) \, dx = \delta_{n,j} .
\end{equation}

The weighted radial Zernike polynomials, denoted by $T_{\text{N,n}}^{\text{p}}(x)$, are defined by the formula
\begin{equation}\label{eq:def:weighted_radial_zernike}
T_{\text{N,n}}^{\text{p}}(x)=x^{\frac{p+1}{2}} \overline{R_{\text{N,n}}^{\text{p}}}(x),
\end{equation}
so that they are orthonormal in the following sense:
\begin{equation}
\int_0^1  {T_{\text{N,n}}^{\text{p}}}(x) {T_{\text{N,j}}^{\text{p}}}(x) \, dx = \delta_{n,j} .
\end{equation}

\subsection{The Gamma Function}

The standard gamma function, denoted by $\Gamma (z)$, is defined by the formula
\begin{equation}\label{eq:def:gammafunc}
  \Gamma(z) = \int_0^\infty x^{z-1} e^{-x} dx.
\end{equation}
For a positive integer $n$,
\begin{equation}
  \Gamma(n+1) = n!.
\end{equation}

\subsection{Bessel Functions}\label{sec:bessel}

The following are well-known properties of Bessel functions of the first kind, denoted by $J_n(x)$. 
The standard definition of Bessel functions can be found, inter alia, in \cite{abramowitz1964handbook}.

\begin{equation}
  e^{\mathrm{i} z \cos{\varphi}} = \sum_{n=-\infty}^\infty \mathrm{i}^n J_n(z) e^{\mathrm{i} n \varphi}
\end{equation}

\begin{equation}
 J_n(z) = z^n \left(\frac{2^{-n}}{\Gamma (n+1)}+O\left(z^2\right)\right)
\end{equation}

\subsection{Surface Harmonics }\label{sec:surfaceharmonics}

In this section we presents properties of integrals on the surface of spheres, that are useful in the polar decomposition of functions
in this paper. A comprehensive discussion of these properties can be found, inter alia, in \cite{slepian1964prolate,serkh2015generalized}.

Let $S_N^m({\bm \xi})$ be the complete set of orthonormal surface harmonics of degree $N=0,1,2,\ldots$, with ${\bm \xi}$ a unit vectors in $D$ dimensions. 
For a hyper-sphere in $D=p+2$ dimensions, there are $h(N,p)$ orthonormal functions in the set, where
\begin{equation}\label{eq:def:hnp}
h(N,p) = (2N+p)\frac{(N+p-1)}{p! N!}
\end{equation}
with $h(N,p)=1$ in the case of $N=p=0$.

\begin{remark}
In the case of $D=2$, the surface harmonics $S_N^m({\bm \xi})$ are the Sine and Cosine functions or the complex exponentials of frequency $N$ (with the proper normalization),
with $h(0,0)=1$ functions for the zero frequency, and $h(N,0)=2$ for the $N>0$ frequency. 
In the case of $D=3$, the surface harmonics are the spherical harmonics $Y_N^m$, with $h(N,0)=2N+1$ functions for the $N$ frequency. 
\end{remark}

Suppose that ${\bm x}= r {\bm \xi}$ and ${\bm y}= r'{\bm n}$, where $r$ and $r'$ are non-negative scalars, and ${\bm \xi}$ and ${\bm n}$ are unit vectors in $D$ dimensions. 
Then,
\begin{equation}\label{eq:exp_bessel}
  \int_\Omega e^{\mathrm{i} c r r' {\bm \xi} \cdot {\bm n}} S_N^m({\bm n}) d\Omega({\bm n})    = H^p_N(crr') S_N^m({\bm \xi})
\end{equation}
where  $ {\bm{x}} \cdot {\bm{t}} = \sum{x_i t_i} $ is the usual inner product,
$\Omega$ is the surface, 
and
\begin{equation}
H^p_N(crr') = \mathrm{i}^N (2 \pi)^{1+p/2} \frac{J_{N+p/2}(crr')}{(crr')^{p/2}},
\end{equation}
where $J_{N+p/2}$ are Bessel functions.
A proof can be found in \cite{slepian1964prolate}.

\section{Analytical Apparatus}\label{sec:analysis}

\subsection{The Truncated Fourier Transform and Associated Integral Operators}\label{sec:truncated_fourier}

In this section we discuss standard definitions of the integral operators associated with GPSFs.
A more detailed discussion of the operators can be found in \cite{slepian1964prolate,serkh2015generalized}.

The unitary Fourier integral transform $\mathcal{F}$ in $D$ dimensions is defined by the formula
\begin{equation}\label{eq:def:Fourier}
  \left(\mathcal{F} f \right)({\bm x}) = \left(\frac{1}{2 \pi }\right)^{D/2} \int_{\mathbb{R}^D} e^{ \mathrm{i}  {\bm{x}} \cdot {\bm{t}}  } f(t)  {d}{\bm t} ,
\end{equation}
where ${\bm t}$ and ${\bm x}$ are vectors, and $ {\bm{x}} \cdot {\bm{t}} = \sum{x_i t_i} $ is the usual inner product.

The truncate Fourier transform $\mathcal{F}_c : L^2(B^D) \rightarrow L^2(c B^D)$ is the restriction of the domain of the Fourier Transform
to functions supported on the unit ball $B^D$ of dimension $D$, centered at the origin, with the range truncated to the ball $c B^D$ of radius $c$,
\begin{equation}\label{eq:def:truncated_Fourier}
  \mathcal{F}_c = {\bm 1}_{c B^D} \mathcal{F} {\bm 1}_{B^D} ,
\end{equation}
where ${\bm 1}_{B^D}$ is the restriction operator. 
In other words,
\begin{equation}
  \left(\mathcal{F}_c f \right)({\bm x}) = \left(\frac{1}{2 \pi }\right)^{D/2} \int_{\|{\bm t}\| \leq 1} e^{ \mathrm{i}  {\bm{x}} \cdot {\bm{t}}  } f(t)  {d}{\bm t} ,
\end{equation}
with ${\bm x}$ restricted to $\|{\bm x}\| \leq 1$.

It is often convenient to scale the truncated Fourier transform to the form ${F}_c : L^2(B^D) \rightarrow L^2(B^D)$ defined by the equation
\begin{equation}\label{eq:scaled_truncated_fourier}
  \left({F}_c \psi \right)(x) = \int_B e^{ \mathrm{i} c \left( {\bm{x}} \cdot {\bm{t}} \right) } \psi(\bm{t})  {d}\bm{t} .
  \end{equation}
In this form, the domain and range are the same function spaces. 
The goal of this paper is to discuss the computation of the GPSFs which are the eigenfunctions of the operator ${F}_c$,
and the associated eigenvalues, i.e. the functions $\psi$ and scalars $\alpha$ that satisfy 
\begin{equation}\label{eq:scaled_truncated_fourier_eig}
 \alpha \psi =  {F}_c \psi .
  \end{equation}

Let $p=D-2$, ${\bm x}= r {\bm \xi}$ and ${\bm y}= r'{\bm n}$, where $r$ and $r'$ are non-negative real numbers, and ${\bm \xi}$ and ${\bm n}$ are unit vectors in $D$ dimensions. Substituting  (\ref{eq:scaled_truncated_fourier}) into (\ref{eq:scaled_truncated_fourier_eig}) with this change of variables yields
\begin{equation}
 \left( {F}_c \psi\right)(r {\bm \xi}) = \alpha \psi(r {\bm \xi})  = \int_0^1 dr' r'^{p+1} \int_\Omega e^{ \mathrm{i} c r r' \left( {\bm{\xi}} \cdot {\bm{n}} \right) } \psi(r' \bm{n})  {d}\Omega(\bm{n}) ,
\end{equation}
and suggests a separation of variables $\psi(r {\bm \xi}) = \Phi(r) S_N^m({\bm \xi})$, such that 
\begin{equation}\label{eq:scaled_truncated_fourier_eig:separation}
 \left( {F}_c \psi\right)(r {\bm \xi}) = \alpha \Phi(r) S_N^m({\bm \xi})  = \int_0^1 dr' r'^{p+1} \Phi(r')  \int_\Omega e^{ \mathrm{i} c r r' \left( {\bm{\xi}} \cdot {\bm{n}} \right) } S_N^m({\bm n}))  {d}\Omega(\bm{n}) .
\end{equation}
Substituting (\ref{eq:exp_bessel}) into (\ref{eq:scaled_truncated_fourier_eig:separation}) yields
\begin{equation}
 \left( {F}_c \psi\right)(r {\bm \xi}) = \alpha \Phi(r) S_N^m({\bm \xi})  = \mathrm{i}^N (2 \pi)^{1+p/2}  \left( \int_0^1 \frac{J_{N+p/2}(crr')}{(crr')^{p/2}} \Phi(r') r'^{p+1} dr' \right) S_N^m({\bm \xi}).
\end{equation}

For a given dimension $D=p+2$ and $c>0$, the radial GPSFs, denoted by $\Phi_{N,n}$, are defined as the eigenfunctions of the radial operator, i.e. the solutions of the integral equation
\begin{equation}\label{eq:def:radial_prol}
  \beta_{N,n} \Phi_{N,n}(r) = \int_0^1 \frac{J_{N+p/2}(crr')}{(crr')^{p/2}} \Phi_{N,n}(r') r'^{p+1} dr' 
\end{equation}
with $\beta_{N,n}$ the eigenvalues of the operator. For each order $N$, the functions are organized so that $\beta_{N,n}>\beta_{N,n+1}$.
Indeed, the GPSFs, denoted by $\psi_{N,n,m}$, are of the form
\begin{equation}\label{eq:def:prol}
  \psi_{N,n,m}(r {\bm \xi}) = \Phi_{N,n}(r) S_N^m({\bm \xi}).
\end{equation}
These functions are the eigenfunctions of $F_c$ (see (\ref{eq:scaled_truncated_fourier_eig})), with the eigenvalues
\begin{equation}\label{eq:alpha}
  \alpha_{N,n,m} = \alpha_{N,n} = \mathrm{i}^N (2 \pi)^{1+p/2} \beta_{N,n},
\end{equation}
and with $N,n = 0,1,\ldots$, and $m=1,2,\ldots,h(N,p)$, where $h(N,p)$ is defined in (\ref{eq:def:hnp}). 

Finally, the weighted radial GPSFs, denoted by $\varphi_{N,n}(r)$, are defined by the formula
\begin{equation}\label{eq:def:weighted_radial_prol}
\varphi_{N,n}(r) = r^{\frac{p+1}{2}} \Phi_{N,n}(r).
\end{equation}
The weighted radial GPSFs satisfy the integral equation
\begin{equation}\label{eq:weighted_radial_prol:eig}
  \gamma_{N,n} \varphi_{N,n}(r) = \int_0^1 {J_{N+p/2}(crr')} \sqrt{crr'} \varphi_{N,n}(r') dr' 
\end{equation}
where 
\begin{equation}\label{eq:gamma}
  \gamma_{N,n} = \beta_{N,n} c^{\frac{p+1}{2}}.
\end{equation}
In other words, the weighted radial GPSFs are the eigenfunctions of the operator $M_{p,c,N}$, defined by the formula
\begin{equation}\label{eq:operator_M}
 \left( M_{p,c,N} \varphi \right)(r) = \int_0^1 {J_{N+p/2}(crr')} \sqrt{crr'} \varphi(r') dr' .
\end{equation}

\subsection{The Associated Differential Operator}

In this section we summarize several properties of the differential operator 
$L_{p,c,N}$, defined by the formula
\begin{equation}\label{eq:operator_L}
   \left( L_{p,c,N} \varphi \right)(x) = \left( \frac{d}{dx}(1-x^2)\frac{d}{dx}\varphi(x) \right) + \left( \frac{1/4-(N+p/2)^2}{x^2} -c^2 x^2\right)\varphi(x).
\end{equation}

In the classic work in \cite{slepian1964prolate}, 
Slepian found that the integral operator $M_{p,c,N}$,  defined in (\ref{eq:operator_M}), commutes with the differential operator $L_{p,c,N}$,
\begin{equation}
  M_{p,c,N} L_{p,c,N} \varphi =   L_{p,c,N} M_{p,c,N} \varphi,
\end{equation}
so that they share the same eigenfunctions (but not eigenvalues).
The eigendecomposition of $L_{p,c,N}$ is therefore 
\begin{equation}\label{eq:operator_L:eig}
   L_{p,c,N} \varphi_{N,n}  = \chi^{p,c,N}_n \varphi_{N,n}, 
\end{equation}
where $\varphi_{N,n}$ are the same functions defined in (\ref{eq:def:weighted_radial_prol}), 
and $\chi^{p,c,N}_n$ the eigenvalues of $L_{p,c,N}$.

\subsection{The Differential Operator in the Basis of Weighted Radial Zernike Polynomials}

The purpose of this section is to discuss the differential operator $L_{p,c,N}$ (defined in (\ref{eq:operator_L})) in the context of the 
weighted radial Zernike polynomials $T_{\text{N,n}}^{\text{p}}$ (defined in (\ref{eq:def:weighted_radial_zernike})). 
The results in this section appear in slightly different forms, or for special cases, in \cite{slepian1964prolate,shkolnisky2007prolate,serkh2015generalized}.

For given $p,c$ and $N$, we define the matrix elements $B_{\text{N}}^{\text{p,c}}(n,n)$ and $B_{\text{N}}^{\text{p,c}}(n-1,n)=B_{\text{N}}^{\text{p,c}}(n,n-1)$
(with $n=0,1,2,\ldots$) of the matrix $B_{\text{N}}^{\text{p,c}}$ by the formulas
\begin{equation}\label{eq:Bmat_nm1n}
\begin{split}
B_{\text{N}}^{\text{p,c}}(n-1,n)&=B_{\text{N}}^{\text{p,c}}(n,n-1)= \\
&=-\frac{c^2 n \left(n+N+\frac{p}{2}\right)}{\sqrt{1-\frac{2}{2 n+N+\frac{p}{2}+1}} \left(2 n+N+\frac{p}{2}\right) \left(2 n+N+\frac{p}{2}+1\right)}, 
\end{split}
\end{equation}
where in the case of $n=0$ we set $B_{\text{N}}^{\text{p,c}}(n-1,n)=B_{\text{N}}^{\text{p,c}}(n,n-1)=0$,
and
\begin{equation}\label{eq:Bmat_nn}
B_{\text{N}}^{\text{p,c}}(n,n)=-\left(\frac{\left(\left(2 n+N+\frac{p}{2}+1\right)
   \left(N+\frac{p}{2}\right)+2 (n+1) n\right) c^2}{\left(2 n+N+\frac{p}{2}\right) \left(2
   n+N+\frac{p}{2}+2\right)}+\kappa _{\text{N,n}}^{\text{p}}\right),
\end{equation}
where in the case of $N=p=n=0$, we define $B_{\text{0}}^{\text{0,c}}(0,0) = -\kappa _{\text{0,0}}^{\text{0}} + c^2/2$,
and with $\kappa _{\text{N,n}}^{\text{p}}$ defined by the formula
\begin{equation}
\kappa _{\text{N,n}}^{\text{p}}=\left(2 n+N+\frac{p}{2}+\frac{1}{2}\right) \left(2
   n+N+\frac{p}{2}+\frac{3}{2}\right) .
\end{equation}

The differential operator $L_{p,c,N}$ (defined in (\ref{eq:operator_L})) applied to $T_{\text{N,n}}^{\text{p}}$ (defined in (\ref{eq:def:weighted_radial_zernike}))
yields a linear combination of $T_{\text{N,n}}^{\text{p}},T_{\text{N,n+1}}^{\text{p}}$ and $T_{\text{N,n-1}}^{\text{p}}$ (or, in the case of $n=0$, a linear combination of 
$T_{\text{N,n}}^{\text{p}}$ and $T_{\text{N,n+1}}^{\text{p}}$), specified in the following equation:
\begin{equation}\label{eq:L_T}
\begin{split}
& \left( L_{p,c,N} T_{\text{N,n}}^{\text{p}} \right)(x) = \\
& = \frac{{d} }{{d} x}\left(\left(1-x^2\right) \frac{{d} T_{\text{N,n}}^{\text{p}}(x)}{{d} x}\right)+\left(\frac{\frac{1}{4}-\left(N+\frac{p}{2}\right)^2}{x^2}-c^2 x^2\right) T_{\text{N,n}}^{\text{p}}(x) = \\
& = B_{\text{N}}^{\text{p,c}}(n,n-1) T_{\text{N,n-1}}^{\text{p}}(x)+B_{\text{N}}^{\text{p,c}}(n,n+1) T_{\text{N,n+1}}^{\text{p}}(x)+B_{\text{N}}^{\text{p,c}}(n,n) T_{\text{N,n}}^{\text{p}}(x)
\end{split}
\end{equation}
An additional analytical verification of this equation is available in the online resources accompanying this paper.

Suppose that ${\bm h}^{p,c,N,n}$ is the vector of coefficients of the weighted radial GPSF $\varphi_{N,n}(r)$ (defined in (\ref{eq:def:weighted_radial_prol})) 
expanded in the basis of weighted radial Zernike polynomials $T_{\text{N,n}}^{\text{p}}$ (defined in (\ref{eq:def:weighted_radial_zernike})), such that 
\begin{equation}\label{eq:weighted_prol_expand}
  \varphi_{N,n}(x) = \sum_{k=0}^{\infty} {h}^{p,c,N,n}_k T_{\text{N,k}}^{\text{p}}(x) ,
\end{equation}
where ${h}^{p,c,N,n}_k$ are the elements of the vector ${\bm h}^{p,c,N,n}$ with $k=0,1,2,\ldots$.
It follows from (\ref{eq:operator_L:eig}), (\ref{eq:L_T}) and (\ref{eq:weighted_prol_expand}) that
${\bm h}^{p,c,N,n}$ are the eigenvectors of the matrix $B_{\text{N}}^{\text{p,c}}$ (defined in (\ref{eq:Bmat_nn}) and (\ref{eq:Bmat_nm1n})),
with the eigenvalues $\chi^{p,c,N}_n$, defined in (\ref{eq:operator_L:eig}):
\begin{equation}\label{eq:operator_L:eig:mat}
\chi^{p,c,N}_n {\bm h}^{p,c,N,n} = B_{\text{N}}^{\text{p,c}} {\bm h}^{p,c,N,n} .
\end{equation}
The eigenvectors are sorted so that $|\chi^{p,c,N}_0| < |\chi^{p,c,N}_1| < \ldots $.

It follows from (\ref{eq:def:weighted_radial_zernike}), (\ref{eq:def:weighted_radial_prol}) and (\ref{eq:weighted_prol_expand}) that the expansion of radial GPSFs $\Phi_{N,n}$ (defined in (\ref{eq:def:radial_prol})), in the basis of normalized radial Zernike polynomials (defined in (\ref{eq:def:norm_radial_zer})) is
\begin{equation}\label{eq:prol_expand}
  \Phi_{N,n}(x) = \sum_{k=0}^{\infty} {h}^{p,c,N,n}_k \overline{R_{\text{N,k}}^{\text{p}}}(x) .
\end{equation}

\subsection{The Relation between $x \frac{{d} \overline{R_{\text{N,n}}^{\text{p}}}(x)}{{d} x}$ and $\overline{R_{\text{N,m}}^{\text{p}}}(x)$ }\label{sec:xdR}

The purpose of this section is to present the relations between $x \frac{{d} \overline{R_{\text{N,n}}^{\text{p}}}(x)}{{d} x}$, 
$\overline{R_{\text{N,m}}^{\text{p}}}(x)$ and $P_{k}^{\left(N+\frac{p}{2},1\right)}\left(1-2 x^2\right)$ which yield an expansion of 
$x \frac{{d} \overline{R_{\text{N,n}}^{\text{p}}}(x)}{{d} x}$ as a linear combination of $\overline{R_{\text{N,m}}^{\text{p}}}(x)$,
with $m=0,1,\ldots,n$.
The following equations are obtained using (\ref{eq:jac1}), (\ref{eq:jac2}), (\ref{eq:jac3}) and (\ref{eq:def:norm_radial_zer}).
An analytic verification of the formulas is available in the on-line resources. 

First, $x \frac{{d} \overline{R_{\text{N,n}}^{\text{p}}}(x)}{{d} x}$ is a linear combination of $\overline{R_{\text{N,n}}^{\text{p}}}(x)$
and $x^N  P_{n-1}^{\left(N+\frac{p}{2},1\right)}\left(1-2 x^2\right)$, as specified by the formula
\begin{equation}
\begin{split}
x &\frac{{d} \overline{R_{\text{N,n}}^{\text{p}}}(x)}{{d} x} = (2 n+N) \overline{R_{\text{N,n}}^{\text{p}}}(x) + \\
 & \left((-1)^{n-1} x^N (2 (n+N)+p) \sqrt{2 \left(2 n+N+\frac{p}{2}+1\right)}\right) P_{n-1}^{\left(N+\frac{p}{2},1\right)}\left(1-2 x^2\right).
\end{split}
\end{equation}

Next, $x^N P_n^{\left(N+\frac{p}{2},1\right)}\left(1-2 x^2\right)$ is a linear combination of $x^N P_{n-1}^{\left(N+\frac{p}{2},1\right)}\left(1-2 x^2\right)$ (with a lower index $n$) and $\overline{R_{\text{N,n}}^{\text{p}}}(x)$, as specified by the formula
\begin{equation}
\begin{split}
&(-1)^n x^N  P_n^{\left(N+\frac{p}{2},1\right)}\left(1-2 x^2\right)= \\
& \frac{(-1)^{n-1} \left(n+N+\frac{p}{2}\right) x^N P_{n-1}^{\left(N+\frac{p}{2},1\right)}\left(1-2 x^2\right)
   +\frac{ \sqrt{2 n+N+\frac{p}{2}+1}}{\sqrt{2}} \overline{R_{\text{N,n}}^{\text{p}}}(x)}{n+N+\frac{p}{2}+1}.
\end{split}
\end{equation}

\subsection{Recurrence Relations between the Eigenvalues of the Integral Operator}

The following formula, due to \cite{serkh2015generalized,greengard2017generalized}, provides the ratio between the eigenvalues $\beta_{N,n}$, defined in (\ref{eq:def:radial_prol})
\begin{equation}\label{eq:eig:rec}
  \frac{\beta_{N,m}}{\beta_{N,n}} = \frac{\int_0^1 x \Phi_{N,n}'(x) \Phi_{N,m}(x) x^{p+1} dx }{\int_0^1 x \Phi_{N,m}'(x) \Phi_{N,n}(x) x^{p+1} dx} .
\end{equation}

\subsection{Relation between Eigenvalues and Expansion Coefficients}

The purpose of this section is to present a relation between the eigenvalues $\gamma_{N,n}$, defined in (\ref{eq:gamma}), and the elements ${h}^{p,c,N,n}_k$, defined in (\ref{eq:weighted_prol_expand}). 
This property is discussed in \cite{slepian1964prolate,serkh2015generalized,greengard2017generalized}.

Substituting (\ref{eq:weighted_prol_expand}) into (\ref{eq:weighted_radial_prol:eig}) yields
\begin{equation}\label{eq:weighted_radial_prol:eig:expand}
  \gamma_{N,n} \sum_{k=0}^{\infty} {h}^{p,c,N,n}_k T_{\text{N,n}}^{\text{p}}(y) = \int_0^1 {J_{N+p/2}(cxy)} \sqrt{cxy} \left(\sum_{k=0}^{\infty} {h}^{p,c,N,n}_k T_{\text{N,n}}^{\text{p}}(x) \right) dx .
\end{equation}

Multiplying the right hand side of (\ref{eq:weighted_radial_prol:eig:expand}) by $y^{-\left(N+\frac{p+1}{2}\right)} \Gamma \left(N+\frac{p}{2}+1\right)$, and considering only one element in the summation over $k$ yields
\begin{equation}\label{eq:expand-eig:rhs}
\begin{split}
& y^{-\left(N+\frac{p+1}{2}\right)} \Gamma \left(N+\frac{p}{2}+1\right) \int_0^1 \sqrt{c x y} T_{\text{N,n}}^{\text{p}}(x) J_{N+\frac{p}{2}}(c x y) \, dx = \\
& = (-1)^n y^{-\left(N+\frac{p+1}{2}\right)} \Gamma \left(N+\frac{p}{2}+1\right) \sqrt{\frac{4 n+2 N+p+2}{c y}} J_{2 n+N+\frac{p}{2}+1}(c y) 
\end{split}
\end{equation}
(see analytic verification in the on-line resources). 
The expansion in a series around $y=0$ yields

\begin{equation}\label{eq:expand-eig:rhs:0}
\begin{split}
&y^{-\left(\frac{p+1}{2}+N\right)} \Gamma \left(N+\frac{p}{2}+1\right) \int_0^1 T_{\text{N,n}}^{\text{p}}(x) J_{N+\frac{p}{2}}(c x y) \sqrt{c x y} \, dx = \\
&= y^{2 n}    \left(\frac{2^{-2 n-N-\frac{p}{2}} c^{2 n+N+\frac{p}{2}+\frac{1}{2}} e^{i n \pi } \Gamma\left(N+\frac{p}{2}+1\right)}{\sqrt{4 n+2 N+p+2} \Gamma \left(2 n+N+\frac{p}{2}+1\right)}+O\left(y^1\right)\right) . 
\end{split}
\end{equation}
In other words, in the limit $y \rightarrow 0^+$, the coefficient ${h}^{p,c,N,n}_k$ controls the behavior of the right hand side of (\ref{eq:weighted_radial_prol:eig:expand}):
\begin{equation}\label{eq:expand-eig:rhs:0:n0}
\begin{split}
 y^{-\left(\frac{p+1}{2}+N\right)} \Gamma \left(N+\frac{p}{2}+1\right) & \int_0^1  T_{\text{N,0}}^{\text{p}}(x) J_{N+\frac{p}{2}}(c x y) \sqrt{c x y} \, dx = \\
& = \frac{2^{-N-\frac{p}{2}}c^{N+\frac{p}{2}+\frac{1}{2}}}{\sqrt{2 N+p+2}}+O\left(y^1\right) .
\end{split}
\end{equation}

Multiplying the left hand side of (\ref{eq:weighted_radial_prol:eig:expand}) by $y^{-\left(N+\frac{p+1}{2}\right)} \Gamma \left(N+\frac{p}{2}+1\right)$, and considering only one element in the summation over $k$ as $y \rightarrow 0^+$ yields
\begin{equation}\label{eq:expand-eig:lhs:0}
\begin{split}
\lim_{y\to 0} \, y^{-\left(N+\frac{p+1}{2}\right)} & \Gamma \left(N+\frac{p}{2}+1\right) T_{\text{N,n}}^{\text{p}}(y) = \\
& \frac{(-1)^n \sqrt{4 n+2 N+p+2} \Gamma \left(n+N+\frac{p}{2}+1\right)}{\Gamma (n+1)}
\end{split}
\end{equation}
(see analytic verification in the on-line resources). 

It follows from (\ref{eq:expand-eig:rhs:0}) and (\ref{eq:expand-eig:lhs:0}) that the eigenvalue is related to the coefficients ${h}^{p,c,N,n}_k$ by the formula
\begin{equation}\label{eq:expand-eig:gam}
  \gamma_{N,n} = \frac{2^{-N-\frac{p}{2}}c^{N+\frac{p}{2}+\frac{1}{2}} {h}^{p,c,N,n}_0 }{ \sqrt{2 N+p+2}  \sum_{k=0}^{\infty} \left( \frac{(-1)^k \sqrt{4 k+2 N+p+2} \Gamma \left(k+N+\frac{p}{2}+1\right)}{\Gamma (k+1)} {h}^{p,c,N,n}_k \right) }.
\end{equation}

\section{Relation to $\mathcal{F}_c$ and Concentration of Energy}\label{sec:native_unscaled}

The eigenfunctions and eigenvalues of ${F}_c$ characterize the operator $\mathcal{F}_c$ through the functions 
$\widetilde{\psi_{N,n,m}}$ and scalars $\nu_{N,n}$ defined in this section; 
together with $\psi_{N,n,m}$, these are the singular functions and singular values of $\mathcal{F}_c$ up to phase).
For a given $p$ and $c$, 
\begin{equation}
 \left( \mathcal{F}_c \psi_{N,n,m}\right) = \nu_{N,n} \widetilde{\psi_{N,n,m}}
\end{equation}
where,
\begin{equation}\label{eq:def:tilde_psi}
  \widetilde{\psi_{N,n,m}}({\bm x}) = c^{-1/2} {\psi_{N,n,m}}({c^{-1}\bm x}),
\end{equation}
and
\begin{equation}\label{eq:nu}
  \nu_{N,n} = \mathrm{i}^N c^{1/2} \gamma_{N,n}  .
\end{equation}
It is convenient to refer to $\psi_{N,n,m}$, $\widetilde{\psi_{N,n,m}}$ and $\nu_{N,n}$ as ``eigenfunctions'' and ``eigenvalues'' of the operator $\mathcal{F}_c$ because of their close relation to the eigenfunctions $\psi_{N,n,m}$ and eigenvalues $\alpha_{N,n},\beta_{N,n}$ and $\gamma_{N,n}$, although the domain of the operator  $\mathcal{F}_c$ is different from its range.

For a given $p,N$, and a growing $c$, the magnitude $|\nu_{N,n}|$ of the first few eigenvalues is very close to $1$. 
In fact, as $c$ grows, there is a growing number of eigenvalues that are numerically indistinguishable from $1$. 
After a certain number of eigenvalues that are very close to $1$, the eigenvalues decay super-algebraically. 
Some examples are presented in the numerical results. 

Consider the Fourier transform $g = \mathcal{F} \psi_{N,n,m}$ of a GPSF,
and the truncated Fourier transform $g_c = \mathcal{F} \psi_{N,n,m}$.
Since the Fourier transform is a unitary operation, and the GPSF is normalized $\|\psi_{N,n,m}\|_2^2=1$,
we also have that 
\begin{equation}
\| \mathcal{F} \psi_{N,n,m} \|_2^2 = \|g\|_2^2=1. 
\end{equation}
Furthermore, since $\psi_{N,n,m}$ is supported on the unit ball by definition, we also have that within the ball $cB^D$ of band $c$, 
\begin{equation}
g({\bm \omega}) = \left( \mathcal{F} \psi_{N,n,m}\right)({\bm \omega}) =  g_c({\bm \omega}) = \left( \mathcal{F}_c \psi_{N,n,m}\right)({\bm \omega}) = \nu_{N,n} \widetilde{\psi_{N,n,m}}({\bm \omega}) ~~,~~ \|{\bm \omega}\| \leq c.
\end{equation}
Therefore, 
\begin{equation}
\| \mathcal{F}_c \psi_{N,n,m} \|_2^2 = \|g_c\|_2^2= |\nu_{N,n}|^2. 
\end{equation}
By definition, the function $g_c({\bm \omega}) = \left(\mathcal{F}_c \psi_{N,n,m}\right)({\bm \omega})$ is identically zero for $\|{\bm \omega}\| \geq c$. 
However, $g({\bm \omega})$ cannot be identically zero on the outside due to classic results in analysis that assert that functions cannot be compactly supported in both the spatial domain and frequency domain. The amount of energy that persists outside the band $c$ is
\begin{equation}
 \| g-g_c \|_2^2 = 1-|\nu_{N,n}|^2.
\end{equation}
It follows that the first few GPSFs, with eigenvalues very close to $1$, are highly concentrated in both the spatial domain 
and the frequency domain. With a small abuse of terminology, we say that these functions are {\em numerically} compactly supported in both the spatial and frequency domains. Subsequent GPSFs with smaller eigenvalues have most of their energy concentrated outside the band $c$. In this sense, the first few GPSFs, with eigenvalues close to $1$, and a rather sharp threshold after which the eigenvalues become very small, are the optimal basis for representing functions that are highly concentrated in both the spatial and frequency domains.

%
%
%
%
%

\section{Algorithms}\label{sec:alg}

\subsection{Computation of GPSFs}

The radial GPSFs, defined in (\ref{eq:def:radial_prol}), are obtained via the expansion (\ref{eq:prol_expand}), in the basis of normalized radial Zernike polynomials (defined in (\ref{eq:def:norm_radial_zer})).
The coefficients of the expansion are the elements of the eigenvectors ${\bm h}^{p,c,N,n}$ of the matrix $B_{\text{N}}^{\text{p,c}}$ defined in (\ref{eq:Bmat_nn}) and (\ref{eq:Bmat_nm1n}). 
The eigenvalues $\chi^{p,c,N}_0$ of the matrix $B_{\text{N}}^{\text{p,c}}$ are sorted so that they grow in magnitude. 
We note that it is not necessary in general to compute all the eigenvectors of the matrix, 
individual eigenvectors can be found analogously to the procedure in \cite{LedermanLaplace2014}.

In the present paper we do not discuss the decay of the elements ${h}^{p,c,N,n}_k$ of the eigenvectors in detail. These elements decay rapidly with $k$, and the matrix $B_{\text{N}}^{\text{p,c}}$ can be truncated without loss of numerical precision. 
In the current implementation we simply ensure that the dimensions of $B_{\text{N}}^{\text{p,c}}$ are sufficiently large to ensure a sufficient number of eigenfunctions and a sufficiently long expansion; 
the appropriate truncation will be discussed in future papers. 

The relation to the ``eigenfunctions'' of the more familiar unscaled Fourier transform
if described in Section \ref{sec:native_unscaled}.
The relation between these radial components of GPSFs and the GPSFs on the D-dimensional ball 
are discussed in Section \ref{sec:truncated_fourier}.

\begin{remark}
Obviously,  the GPSFs computed through the eigendecomposition above have a degree of freedom in their sign (if $\Phi_{N,n}(x)$ is a normalized eigenfunction, then $-\Phi_{N,n}(x)$ is also a normalized eigenfunction). 
For the sake of consistency across different implementations, we remove the ambiguity by setting the sign of the first element ${h}^{p,c,N,n}_0$ of the eigenvectors to be positive for even $n$ and negative for odd $n$. 
\end{remark}

\begin{remark}\label{remark:eig_alg}
The sign standardization above, and the accurate computation of eigenvalues using  (\ref{eq:gamma}) in the next section rely in some cases on obtaining relative precision in some elements of the eigenvectors ${\bm h}^{p,c,N,n}$. 
While relative precision can be obtained in the cases in question (see \cite{osipov2017evaluation}), not all eigedecomposition algorithms achieve it. 
\end{remark}

\subsection{Computation of the First Associated Eigenvalue}

The first eigenvalue $\gamma_{N,0}$, defined in (\ref{eq:gamma}), is computed using the relation (\ref{eq:expand-eig:gam}), with $n=0$.
Alternatively, it can often be computed via (\ref{eq:weighted_radial_prol:eig}) by  numerical integration.

\subsection{Computation of Associated Eigenvalues}

The expansion $\widetilde{{\bm h}^{p,c,N,n}}$ of $x \Phi_{N,n}'(x)$ such that 
\begin{equation}
  x \Phi_{N,n}'(x) = \sum_{k=0}^{\infty} \widetilde{{h}^{p,c,N,n}_k} \overline{R_{\text{N,k}}^{\text{p}}}(x) ,
\end{equation}
is computed from the eigenvectors ${{\bm h}^{p,c,N,n}}$ recovered in previous steps of the algorithm,  using the equations in Section \ref{sec:xdR}.

It follows from (\ref{eq:eig:rec}) that the ratio between consecutive eigenvalues is related to the elements  ${{\bm h}^{p,c,N,n}}$ of the expansion of  $\Phi_{N,n}'(x)$ and the elements $\widetilde{{\bm h}^{p,c,N,n}}$ of the expansion of $x \Phi_{N,n}'(x)$ via the following formula
\begin{equation}
  \frac{\gamma_{N,n+1}}{\gamma_{N,n}} = \frac{ \widetilde{{\bm h}^{p,c,N,n}}\cdot {{\bm h}^{p,c,N,n+1}}  }{  \widetilde{{\bm h}^{p,c,N,n+1}}\cdot {{\bm h}^{p,c,N,n}} } ,
\end{equation}
where $\cdot$ is the standard inner product between vectors.
This ratio between consecutive eigenvalues is used to compute the subsequent eigenvalues once $\gamma_{N,0}$ is computed by the means described in the previous section.

The other associated eigenvalues  $\alpha_{N,n,m} = \alpha_{N,n}$, $\beta_{N,n}$ and 
$\nu_{N,n}$ are computed from $\gamma_{N,n}$ using equations (\ref{eq:alpha}), (\ref{eq:def:radial_prol}) and (\ref{eq:nu}).

\section{Numerical Results}\label{sec:results}

The algorithms describe above have been implemented in MATLAB\textsuperscript{TM}.
The following figures are examples of GPSFs and associated eigenvalues. 
The code for reproducing these plots is available at \url{http://github.com/lederman/prol}.

\begin{remark}
For the sake of simplicity, the current MATLAB\textsuperscript{TM} implementation makes use of MATLAB's ``eig'' eigendecomposition.
This subroutine retains relative precision only in certain elements, but not in all elements, as discussed in Remark \ref{remark:eig_alg}.
The truncation is corrected for certain computations using an inverse power method step on the first eigenvector. 
For this reason and other reasons, future implementations will replace the eigendecomposition procedure. 
\end{remark}

\begin{figure}[h]
\includegraphics[width=\textwidth]{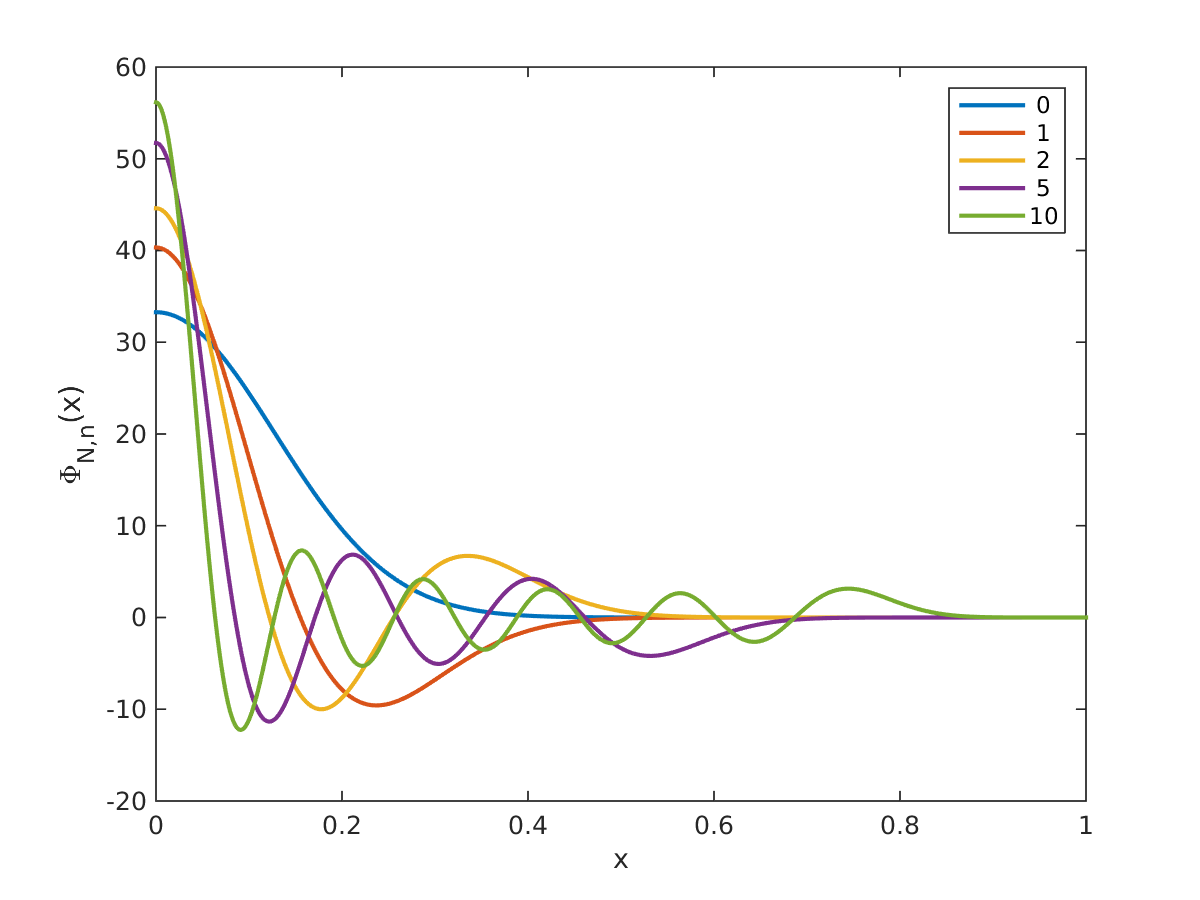}
\caption{Examples of radial GPSFs, with $D=3$, $c=20\pi$, and $N=0$.  }
\end{figure}

\begin{figure}[h]
\includegraphics[width=\textwidth]{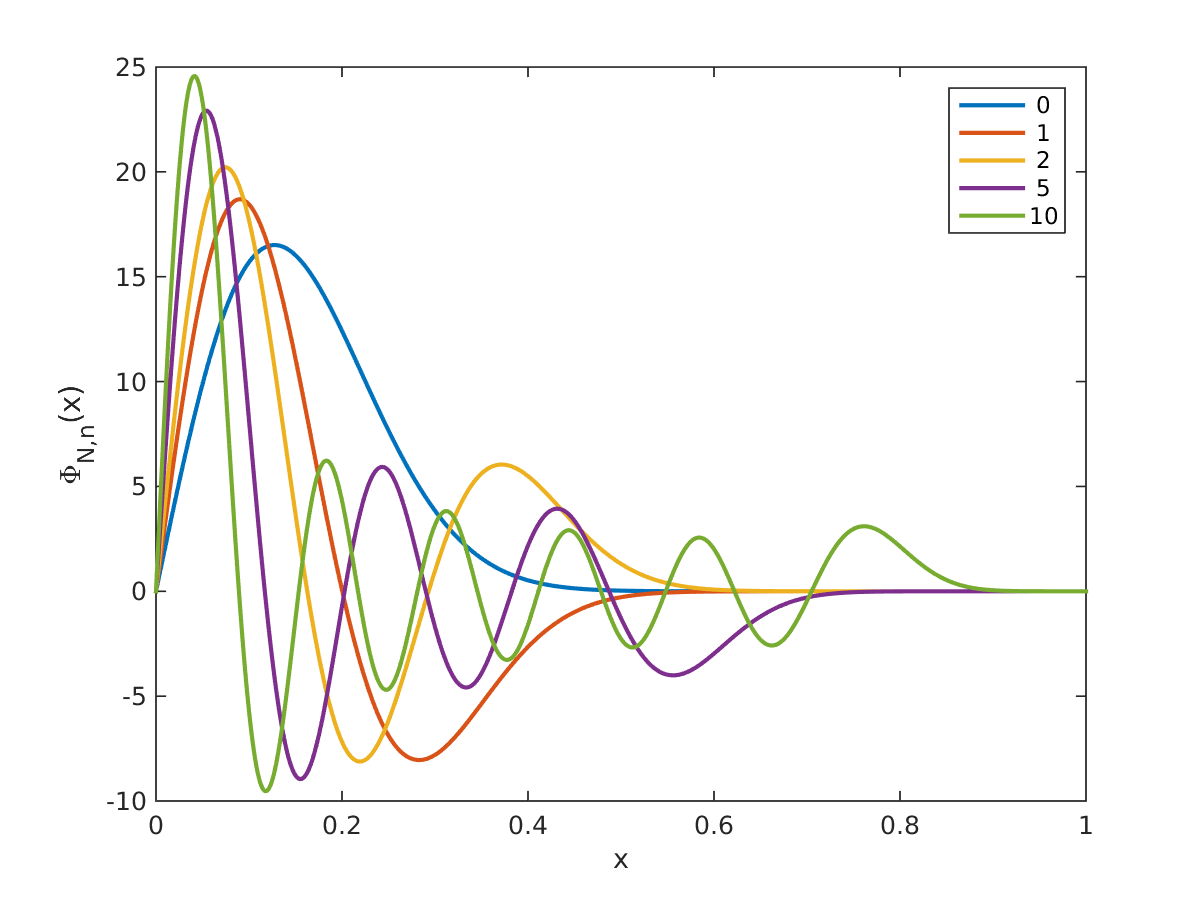}
\caption{Examples of radial GPSFs, with $D=3$, $c=20\pi$ and $N=1$.  }
\end{figure}

\begin{figure}[h]
\includegraphics[width=\textwidth]{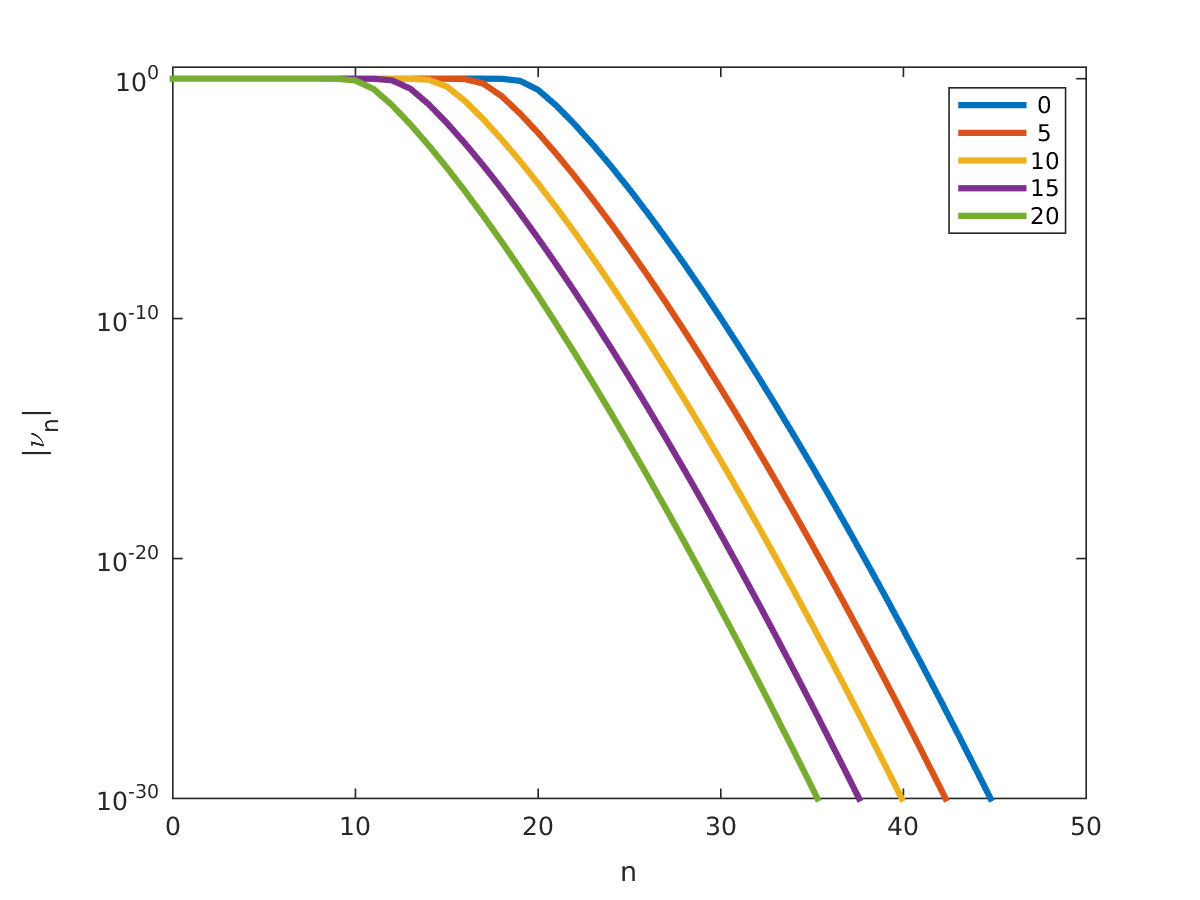}
\caption{Magnitude of eigenvalues $|\nu_{N,n}|$, with $D=3$ and $c=20\pi$, for different values of $N$  }
\end{figure}

\begin{figure}[h]
\includegraphics[width=\textwidth]{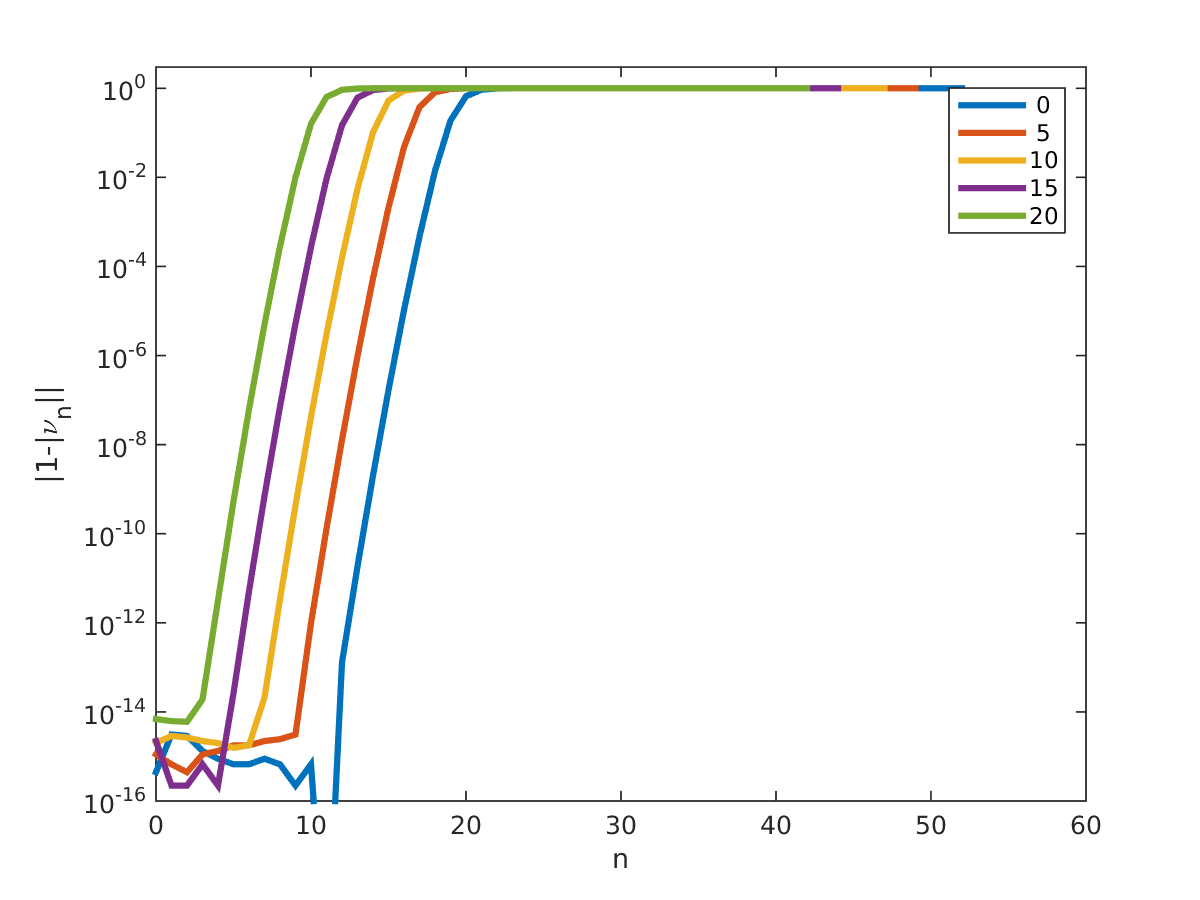}
\caption{$1-|\nu_{N,n}|$, with $D=3$ and $c=20\pi$, for different values of $N$, demonstrating that the magnitude of the first eigenvalues is numerically indistinguishable from one. The difference from one in the first few eigenvalues is due to numerical precision.  }
\end{figure}

\clearpage

\section{Conclusions}\label{sec:conclusions}

In this paper we reviewed the essential results for numerical computation of the eigenfunctions of the 
truncated Fourier transform, which is restricted to the unit ball in the spatial domain and band-limited in the frequency domain. 
The paper is accompanied by an example implementation, available as an open-source project at \url{http://github.com/lederman/prol}.

The code and paper are preliminary results of the ongoing project in \url{http://github.com/lederman/prol}. 
A more comprehensive code in additional languages is being developed. 
We welcome comments and contributions to both the code and the text through the open-source project.

\section{Acknowledgements}

We would like to thank Philip Greengard, Boris Landa and Kirill Serkh for their help. 
The author was partially supported by Award Number R01GM090200 from the NIGMS to Amit Singer.

\bibliography{bib}{}

\begin{thebibliography}{10}

\bibitem{abramowitz1964handbook}
Milton Abramowitz and Irene~A Stegun.
\newblock {\em Handbook of mathematical functions: with formulas, graphs, and
  mathematical tables}.
\newblock National Bureau of Standards, New York, 1966.

\bibitem{bertero1985commuting}
Mario Bertero and F.~Alberto Grunbaum.
\newblock Commuting differential operators for the finite {L}aplace transform.
\newblock {\em Inverse Problems}, 1(3):181, 1985.

\bibitem{greengard2017generalized}
Philip Greengard and Kirill Serkh.
\newblock On generalized prolate spheroidal functions.
\newblock {\em Manuscript in preparation}, 2017.

\bibitem{landa2017approximation}
Boris Landa and Yoel Shkolnisky.
\newblock Approximation scheme for essentially bandlimited and
  space-concentrated functions on a disk.
\newblock {\em Applied and Computational Harmonic Analysis}, 43(3):381--403,
  2017.

\bibitem{landa2017steerable}
Boris Landa and Yoel Shkolnisky.
\newblock Steerable principal components for space-frequency localized images.
\newblock {\em SIAM Journal on Imaging Sciences}, 10(2):508--534, 2017.

\bibitem{LedermanLaplace2014}
Roy~R. Lederman.
\newblock {\em On the Analytical and Numerical Properties of the Truncated
  {L}aplace Transform}.
\newblock PhD thesis, Yale University, 2014.

\bibitem{LedermanLaplace2}
Roy~R. Lederman and Vladimir Rokhlin.
\newblock On the analytical and numerical properties of the truncated {L}aplace
  transform {I}.
\newblock {\em SIAM Journal on Numerical Analysis}, 53(3):1214--1235, 2015.

\bibitem{lederman2016representation}
Roy~R. Lederman and Amit Singer.
\newblock A representation theory perspective on simultaneous alignment and
  classification.
\newblock {\em arXiv preprint arXiv:1607.03464}, 2016.

\bibitem{lederman2017continuously}
Roy~R. Lederman and Amit Singer.
\newblock Continuously heterogeneous hyper-objects in cryo-em and 3-d movies of
  many temporal dimensions.
\newblock {\em arXiv preprint arXiv:1704.02899}, 2017.

\bibitem{lederman2017lower}
Roy~R. Lederman and Stefan Steinerberger.
\newblock Lower bounds for truncated {F}ourier and {L}aplace transforms.
\newblock {\em Integral Equations and Operator Theory}, 87(4):529--543, 2017.

\bibitem{osipov2017evaluation}
Andrei Osipov.
\newblock Evaluation of small elements of the eigenvectors of certain symmetric
  tridiagonal matrices with high relative accuracy.
\newblock {\em Applied and Computational Harmonic Analysis}, 43(2):173--211,
  2017.

\bibitem{osipov2013prolate}
Andrei Osipov, Vladimir Rokhlin, and Hong Xiao.
\newblock Prolate spheroidal wave functions of order zero.
\newblock {\em Springer Ser. Appl. Math. Sci}, 187, 2013.

\bibitem{rokhlin2007approximate}
Vladimir Rokhlin and Hong Xiao.
\newblock Approximate formulae for certain prolate spheroidal wave functions
  valid for large values of both order and band-limit.
\newblock {\em Applied and Computational Harmonic Analysis}, 22(1):105--123,
  2007.

\bibitem{serkh2015generalized}
Kirill Serkh.
\newblock On generalized prolate spheroidal functions.
\newblock Technical report, Technical Report TR-1519, Department of
  Mathematics, Yale University, 2015.

\bibitem{shkolnisky2007prolate}
Yoel Shkolnisky.
\newblock Prolate spheroidal wave functions on a disc—integration and
  approximation of two-dimensional bandlimited functions.
\newblock {\em Applied and Computational Harmonic Analysis}, 22(2):235--256,
  2007.

\bibitem{shkolnisky2006approximation}
Yoel Shkolnisky, Mark Tygert, and Vladimir Rokhlin.
\newblock Approximation of bandlimited functions.
\newblock {\em Applied and Computational Harmonic Analysis}, 21(3):413--420,
  2006.

\bibitem{slepian1964prolate}
David Slepian.
\newblock Prolate spheroidal wave functions, {F}ourier analysis and uncertainty
  {- IV}: extensions to many dimensions; generalized prolate spheroidal
  functions.
\newblock {\em Bell Labs Technical Journal}, 43(6):3009--3057, 1964.

\bibitem{slepian1961prolate}
David Slepian and Henry~O Pollak.
\newblock Prolate spheroidal wave functions, {F}ourier analysis and uncertainty
  {- I}.
\newblock {\em Bell Labs Technical Journal}, 40(1):43--63, 1961.

\bibitem{xiao2003high}
Hong Xiao and Vladimir Rokhlin.
\newblock High-frequency asymptotic expansions for certain prolate spheroidal
  wave functions.
\newblock {\em Journal of Fourier Analysis and Applications}, 9(6):575--596,
  2003.

\bibitem{xiao2001prolate}
Hong Xiao, Vladimir Rokhlin, and Norman Yarvin.
\newblock Prolate spheroidal wavefunctions, quadrature and interpolation.
\newblock {\em Inverse problems}, 17(4):805, 2001.

\end{thebibliography}
\bibliographystyle{plain}

\end{document}